\documentclass[12pt,a4paper,reqno]{amsart}
\usepackage{amssymb}
\usepackage{dcpic,pictex}
\usepackage{multirow}
\usepackage{graphicx}
\usepackage{slashbox}
\usepackage{cite}

\usepackage{hyperref}
\usepackage[main=english,french]{babel}
\newtheorem{definition}{Definition}[section]
\newtheorem{example}{Example}[section]
\newtheorem{theorem}{Theorem}[section]

\newcommand{\D}{\mathbb{D}}

\newcommand{\I}{\mathbb{I}}
\newcommand{\R}{\mathbb{R}}
\newcommand{\C}{\mathbb{C}}
\newcommand{\Com}{\mathbb{C}}
\newcommand{\N}{\mathbb{N}}

\begin{document}


\title[Convolution series]{Convolution series and the generalized convolution Taylor formula}

\author{Yuri Luchko}
\curraddr{Beuth Technical University of Applied Sciences Berlin,  
     Department of  Mathematics, Physics, and Chemistry,  
     Luxemburger Str. 10,  
     13353 Berlin,
Germany}
\email{luchko@beuth-hochschule.de}

\subjclass[2010]{26A33; 26B30; 44A10; 45E10}
\dedicatory{}
\keywords{Convolution series, convolution polynomials, Sonine kernel, general fractional derivative, general fractional integral, sequential general fractional derivative, generalized convolution Taylor formula}

\begin{abstract}
In this paper, we deal with the convolution series that are a far reaching generalization of the conventional power series and the power series with the fractional exponents including the Mittag-Leffler type functions. Special attention is given to the most interesting case of the convolution series generated by the Sonine kernels. In this paper, we first formulate and prove the second fundamental theorem for the general fractional integrals and the $n$-fold general sequential  fractional derivatives of both the Riemann-Liouville and the Caputo types. These results are then employed for derivation of two different forms of a generalized convolution Taylor formula for representation of a function as a convolution polynomial with a remainder in form of a composition of the $n$-fold general fractional integral and the $n$-fold general sequential fractional derivative of the Riemann-Liouville and the Caputo types, respectively. We also discuss the generalized Taylor series in form of convolution series and deduce the formulas for its coefficients in terms of the $n$-fold general sequential  fractional derivatives. 
\end{abstract}

\maketitle


\section{Introduction}\label{sec:1}

The power series are a very important instrument both in mathematics and its applications. Without any loss of generality, a power series can be represented in the form
\begin{equation}
\label{I-0}
\Sigma(x) = \sum_{j=0}^{+\infty} a_j\, h_{j+1}(x),\ a_j \in \R,
\end{equation}
where with $h_\beta$ we denote the following power function 
\begin{equation}
\label{I-2}
h_{\beta}(x)=\frac{x^{\beta-1}}{\Gamma(\beta)},\ \ \beta >0.
\end{equation}
A given function $f$ ($f\in C^{n-1}(-\delta,\delta)$, $\delta>0$) can be approximated by its Taylor polynomial through the well-known Taylor formula 
\begin{equation}
\label{I-1}
f(x) = \sum_{j=0}^{n-1} a_j\, h_{j+1}(x)\ + \ r_n(x),\ \ a_j = f^{(j)}(0),\ j=0,1,\dots,n-1 
\end{equation}
with the remainder $r_n=r_n(x)$ in different forms including the one suggested by Lagrange (for the functions $f\in C^{n}(-\delta,\delta)$):
 \begin{equation}
\label{I-3}
r_n(x) = f^{(n)}(\xi)\, h_{n+1}(x),\ \ \xi\in (0,x). 
\end{equation}
In the case of a function $f\in C^{\infty}(-\delta,\delta)$, its Taylor series is introduced by letting $n$ go to $+\infty$ in the formula \eqref{I-1}:
\begin{equation}
\label{I-4}
f(x) = \sum_{j=0}^{+\infty} a_j\, h_{j+1}(x),\ \ a_j = f^{(j)}(0).
\end{equation}

To illustrate the main objects of this paper, let us represent the power function 
$h_{j+1}$ from the formulas \eqref{I-0}, \eqref{I-1}, and \eqref{I-4} as a convolution power. We start with the  
following  important formula that is a direct consequence from the well-known representation of the Euler Beta-function in terms of the Gamma-function:
\begin{equation}
\label{2-9}
(h_{\alpha} \, * \, h_\beta)(x) \, = \, h_{\alpha+\beta}(x),\ \alpha,\beta >0,\ x>0,
\end{equation}
where $*$ is the  Laplace convolution
\begin{equation}
\label{2-2}
(f\, *\, g)(x) = \int_0^{x}\, f(x-\xi)g(\xi)\, d\xi.
\end{equation}
Then,   by      repeatedly applying   the formula \eqref{2-9}, we get the representation
\begin{equation}
\label{2-13}
h_{\alpha}^{<n>}(x) = h_{n\alpha}(x),\ n\in \N,
\end{equation}
where the expression $g^{<n>}$ stands for a convolution power
\begin{equation}
\label{I-6}
g^{<n>}(x):= \begin{cases}
1,& n=0,\\
g(x), & n=1,\\
(\underbrace{g* \ldots * g}_{n\ \mbox{times}})(x),& n=2,3,\dots .
\end{cases}
\end{equation}
Setting $\alpha =1$ in the formula \eqref{2-13} leads to the representation
\begin{equation}
\label{I-5}
h_{1}^{<n>}(x) = h_{n}(x),\ j=0,1,2\dots.
\end{equation}
Denoting $h_{1}$  with $\kappa$, the formulas \eqref{I-0}, \eqref{I-1}, and \eqref{I-4} can be rewritten as follows:
\begin{equation}
\label{I-0-1}
\Sigma_\kappa(x) = \sum_{j=0}^{+\infty} a_j\, \kappa^{<j+1>}(x),
\end{equation}
\begin{equation}
\label{I-1-1}
f(x) = \sum_{j=0}^{n-1} a_j\, \kappa^{<j+1>}(x)\ + \ r_n(x),\ \ a_j = a(f,\kappa,j), 
\end{equation}
\begin{equation}
\label{I-4-1}
f(x) = \sum_{j=0}^{+\infty} a_j\, \kappa^{<j+1>}(x),\ \ a_j = a(f,\kappa,j).
\end{equation}

The series in form \eqref{I-0-1} with the functions $\kappa$  that are continuous on the positive real semi-axis and can have an integrable singularity at the origin were recently introduced and studied in \cite{Luc21b} in the framework of an operational calculus for the general fractional derivatives with the Sonine kernels. In \cite{Luc21b}, these series were called convolution series. In the case, the function $\kappa$ is a Sonine kernel, solutions to the initial-value problems for the linear single- and multi-term fractional differential equations with the $n$-fold general fractional derivatives were derived in \cite{Luc21b} in terms of the convolution series \eqref{I-0-1}.

The main focus of this paper is on the generalized convolution Taylor formula in form \eqref{I-1-1} and the generalized convolution Taylor series in form \eqref{I-4-1} with the Sonine kernels $\kappa$. It is worth mentioning that some generalizations of the Taylor formula \eqref{I-1} and the Taylor series \eqref{I-4} to the case of the power series with the fractional exponents have been already considered in the literature (see \cite{Tru} and the references therein). In this paper, we obtain some formulas of this kind as a particular case of our general results while specifying them for the Sonine kernel $\kappa(x) = h_\alpha(x),\ 0<\alpha <1$. 

The rest of the paper is organized as follows. In the 2nd Section, we introduce the main tools needed for derivation of our results including the suitable spaces of functions, the Sonine kernels, and the general fractional integrals and derivatives with the Sonine kernels. The 3rd Section is devoted to the generalized Taylor series in form of  a convolution series generated by the Sonine kernels. In particular, the formulas for the coefficients of the generalized convolution Taylor series in terms of the generalized fractional integrals and derivatives are derived and some particular cases are discussed in details. In the 4th Section, the convolution Taylor polynomials, approximations of the functions from some special spaces in terms of these polynomials (generalized convolution Taylor formula), as well as the formulas for a remainder in the generalized convolution Taylor formula are presented.  

\section{Preliminary results}

Recently, the general fractional integrals (GFIs) and the general fractional derivatives (GFDs) with the Sonine kernels became a subject of active research in Fractional Calculus (FC), see \cite{Han20,Koch11,Koch19_1,LucYam16,LucYam20,Luc21a,Luc21b,Luc21c,Tar}. In this paper, we first derive the second fundamental theorem for the $n$-fold general sequential  fractional derivatives of both the Riemann-Liouville and the Caputo types and then employ them for  derivation of the generalized convolution Taylor formulas in two different forms. 

The GFD in the Riemann-Liouville sense is defined in form of the following integro-differential operator of convolution type
\begin{equation}
\label{FDR-L}
(\D_{(k)}\, f) (x) = \frac{d}{dx}(k\, *\, f)(x) = \frac{d}{dx} \int_0^t k(x-\xi)f(\xi)\, d\xi,
\end{equation}
whereas in the definition of the GFD in the Caputo sense the order of differentiation and integration is interchanged:
\begin{equation}
\label{FD-C}
(\,_*\D_{(k)}\, f) (x) = (k\, *\, f^\prime)(x) =  \int_0^t k(x-\xi)f^\prime(\xi)\, d\xi.
\end{equation}

The prominent particular cases of the GFDs \eqref{FDR-L} and \eqref{FD-C} are the Riemann-Liouville  and the Caputo fractional derivatives of order $\alpha,\ 0< \alpha < 1$ ($k(x)= h_{1-\alpha}(x), \ 0< \alpha < 1$), the multi-term Riemann-Liouville  and Caputo fractional derivatives ($k(x) = \sum_{k=1}^n a_k\, h_{1-\alpha_k}(x),
\ \ 0 < \alpha_1 <\dots < \alpha_n < 1,\ a_k \in \R,\ k=1,\dots,n$), and the Riemann-Liouville  and Caputo fractional derivatives of distributed order ($k(x) = \int_0^1 h_{1-\alpha}(x)\, d\rho(\alpha)$ with a Borel measure $\rho$ defined on the interval $[0,\, 1]$). For other particular cases see \cite{LucYam20,Luc21a} and the references therein. For applications of the GFDs in general fractional dynamics and in general non-Markovian quantum dynamics we refer the interested readers to the recent publications \cite{Tar1,Tar2}. 

It is worth mentioning that the Laplace convolution of the kernel $\kappa(x) = h_{\alpha}(x), \ \alpha >0$ with a function $f$ is known as the Riemann-Liouville  fractional integral of order $\alpha$:
\begin{equation}
\label{RLI}
\left(I^\alpha_{0+} f \right) (x)=(h_\alpha \, * f)(x) = \frac{1}{\Gamma (\alpha )}\int\limits_0^x(x -\xi)^{\alpha -1}f(\xi)\,d\xi,\ x>0.
\end{equation}
Thus, the   Riemann-Liouville  fractional derivative of order $\alpha,\ 0< \alpha < 1$ can be represented as a composition of the first-order derivative and the Riemann-Liouville fractional integral of order $1-\alpha$:
\begin{equation}
\label{RLD}
(D_{0+}^\alpha\, f)(x) = \frac{d}{dx}(h_{1-\alpha} \, * f)(x)=\frac{d}{dx}(I^{1-\alpha}_{0+} f)(x),\ x>0.
\end{equation}
The Caputo fractional derivative is the Riemann-Liouville fractional integral of order $1-\alpha$ applied to the first order derivative of a function $f$:
\begin{equation}
\label{CD}
(\,_*D_{0+}^\alpha\, f)(x) = (h_{1-\alpha} \, * f^\prime)(x)=(I^{1-\alpha}_{0+} f^\prime)(x),\ x>0.
\end{equation}



Applying the formula \eqref{2-9} to the kernels $h_{\alpha}(x)$ and $h_{1-\alpha}(x)$ of the Riemann-Liouville  fractional integral and the Riemann-Liouville  or Caputo fractional derivatives, respectively, we get the relation 
\begin{equation}
\label{3-1}
(h_{\alpha} \, *\, h_{1-\alpha} )(x) = h_1(x)  = \{1\},\ 0<\alpha <1,\ x>0,
\end{equation}
where by $ \{1\}$ we denote the function that is identically equal to $1$. 

In \cite{Son}, the general kernels $\kappa,\, k$ that satisfy the  condition 
\begin{equation}
\label{3-2}
(\kappa \, *\, k )(x) = h_1(x)  = \{1\},\ x>0
\end{equation}
have been introduced and studied.

Nowadays, the functions that satisfy the condition \eqref{3-2} (Sonine condition) are called the Sonine kernels. For a Sonine kernel $\kappa$, the kernel $k$ that satisfies the Sonine condition \eqref{3-2} is called its associated kernel. Of course, $\kappa$ is then an associated kernel to $k$. The set of the Sonine kernels is denoted by $\mathcal{S}$.

It was known already to Sonine (\cite{Son}) that the GFD \eqref{FDR-L} with the Sonine kernel $k$ and the GFI with its associated kernel $\kappa$ in form
\begin{equation}
\label{GFI}
(\I_{(\kappa)}\, f) (x) =  (\kappa\, *\, f)(x) = \int_0^t \kappa(x-\xi)f(\xi)\, d\xi
\end{equation}
satisfy the so called 1st fundamental theorem of FC (see \cite{Luc20}), i.e., 
\begin{equation}
\label{FTL_S}
(\D_{(k)}\, \I_{(\kappa)}\, f) (x) = \frac{d}{dx}(k\, *\, (\kappa \, *\, f))(x) = 
\frac{d}{dx}(\{1\}\, *\  f)(x) = f(x)
\end{equation}
on the corresponding spaces of functions. 

In the publications \cite{Luc21b,Han20,Koch11,Luc21a,Luc21c}, some important classes of the Sonine kernels as well as  the GFDs and GFIs with these kernels on the appropriate spaces of functions have been introduced and studied. 




In \cite{Han20}, the case of the Sonine kernels in form of the singular (unbounded in a neighborhood of the point zero) locally integrable
completely monotone functions was discussed. A typical example of a pair of the Sonine kernels of this sort is as follows (\cite{Han20}):
\begin{equation}
\label{3-8} 
\kappa(x) = h_{1-\beta+\alpha}(x)\, +\, h_{1-\beta}(x),\ 0<\alpha < \beta <1,
\end{equation}
\begin{equation}
\label{3-9} 
k(x) = x^{\beta -1}\, E_{\alpha,\beta}(-x^\alpha),
\end{equation}
where $E_{\alpha,\beta}$ stands for the two-parameters Mittag--Leffler function that is defined by the following convergent series:
\begin{equation}
\label{ML}
E_{\alpha,\beta}(z) = \sum_{k=0}^{+\infty} \frac{z^k}{\Gamma(\alpha\, k + \beta)},\ \alpha >0,\ \beta,z\in \Com.
\end{equation}

In the publications \cite{Luc21b,Luc21a,Luc21c}, the Sonine kernels continuous on $\R_+$ that have an integrable singularity at the origin were treated. In this paper, we mainly deal with the kernels from this space of functions that is defined as follows:
\begin{equation}
\label{subspace}
 C_{-1,0}(0,+\infty) \, = \, \{f:\ f(x) = x^{p-1}f_1(x),\ x>0,\ 0<p<1,\ f_1\in C[0,+\infty)\}.
\end{equation}
\begin{definition}[\cite{Luc21c}]
\label{dd2}
Let $\kappa,\, k \in C_{-1,0}(0,+\infty)$ be a pair of the Sonine kernels, i.e., the Sonine condition \eqref{3-2} be fulfilled. The set of such Sonine kernels    is denoted by $\mathcal{L}_{1}$:
\begin{equation}
\label{Son}
(\kappa,\, k \in \mathcal{L}_{1} ) \ \Leftrightarrow \ (\kappa,\, k \in C_{-1,0}(0,+\infty))\wedge ((\kappa\, *\, k)(x) \, = \, h_1(x)).
\end{equation}
\end{definition}



For derivation of our results, we employ some properties of the GFI \eqref{GFI} and the GFDs \eqref{FDR-L} and \eqref{FD-C} of the Riemann-Liouville and the Caputo type, respectively, with the kernels $\kappa$ and $k$ from the class $\mathcal{L}_{1}$ on the space of functions $C_{-1}(0,+\infty)$ and its subspaces. In the rest of this section, we present some formulas derived in \cite{Luc21b,Luc21a,Luc21c} that are needed for the further discussions.

First, we mention  the relations 
\begin{equation}
\label{GFI_11}
(I^0_{0+}\, f)(x) = f(x),\ \ (I^1_{0+}\, f)(x) = \int_0^x\, f(\xi)\, d\xi 
\end{equation} 
valid for the Riemann-Liouville fractional integral \eqref{RLI}. Combining the formulas \eqref{GFI} and \eqref{GFI_11}, it is natural to define the GFIs with the kernels $\kappa(x) = h_0(x)$ and $\kappa(x) = h_1(x)$ (that do not belong to the set $\mathcal{L}_{1}$) as follows:
\begin{equation}
\label{GFI_12}
(\I_{(h_0)}\, f)(x) := (I^0_{0+}\, f)(x) = f(x),\ \ (\I_{(h_1)}\, f)(x) := (I^1_{0+}\, f)(x) = \int_0^x\, f(\xi)\, d\xi. 
\end{equation} 

In the formula \eqref{GFI_12}, the function $h_0$ is interpreted as a kind of the $\delta$-function that plays the role of a unity with respect to multiplication in form of the Laplace convolution. In particular, one can extend  the formula \eqref{3-1} valid for $0<\alpha<1$ in the usual sense to the case $\alpha =1$ that has to be understand in the sense of generalized functions:
\begin{equation}
\label{3-1-1}
(h_{1}\, * \, h_{0})(x) := h_1(x),\ t>0.
\end{equation}

In its turn, the convolution formula \eqref{3-1-1} suggests the following natural definitions of the GFDs with the kernels $k = h_1(t)$ 
and $k(t) = h_0(t)$, respectively:
\begin{equation}
\label{GFD_12_1}
( \D_{(h_1)}\, f)(x):= \frac{d}{dx} ( \I_{(h_0)}\, f)(x) = \frac{d}{dx} ( I_{0+}^0\, f)(x) = f^\prime(x), 
\end{equation} 
\begin{equation}
\label{GFD_12_2}
( \D_{(h_0)}\, f)(x):= \frac{d}{dx} ( \I_{(h_1)}\, f)(x) = \frac{d}{dx} ( I_{0+}^1\, f)(x) = f(x), 
\end{equation}  
\begin{equation}
\label{GFD_12_3}
(\,_* \D_{(h_1)}\, f)(x):= ( \I_{(h_0)}\, f^\prime)(x) =  ( I_{0+}^0\, f^\prime)(x) = f^\prime(x), 
\end{equation} 
\begin{equation}
\label{GFD_12_4}
(\,_* \D_{(h_0)}\, f)(x):= ( \I_{(h_1)}\, f^\prime)(x) = ( I_{0+}^1\, f^\prime)(x) = f(x)-f(0). 
\end{equation} 

Several other particular cases of the GFI \eqref{GFI} can be easily constructed using the known Sonine kernels (see \cite{Luc21a}). Here,  we mention just one of them:
\begin{equation}
\label{GFI_2}
(\I_{(\kappa)}\, f)(x) = (I^{1-\beta+\alpha}_{0+}\, f)(x) + (I^{1-\beta}_{0+}\, f)(x) ,\ x>0
\end{equation} 
with the Sonine kernel $\kappa$ defined by \eqref{3-8}.

The GFD of the Riemann-Liouville type that corresponds to the GFI \eqref{GFI_2} has the Mittag-Leffler function \eqref{3-9} in the kernel:
\begin{equation}
\label{GFD_2}
(\D_{(k)}\, f)(x) = \frac{d}{dx}\, \int_0^x (x-\xi)^{\beta -1}\, E_{\alpha,\beta}(-(x-\xi)^\alpha)\, f(\xi)\, d\xi,\ 0<\alpha < \beta <1,\ x>0.
\end{equation} 
The GFD of the Caputo type with the Mittag-Leffler function \eqref{3-9} in the kernel takes the form
\begin{equation}
\label{GFD_3}
(\,_*\D_{(k)}\, f)(x) =  \int_0^x (x-\xi)^{\beta -1}\, E_{\alpha,\beta}(-(x-\xi)^\alpha)\, f^\prime(\xi)\, d\xi,\ 0<\alpha < \beta <1,\ x>0.
\end{equation} 

Several important properties of the GFI \eqref{GFI} on the space $C_{-1}(0,+\infty)$ can be easily derived from the known properties of the Laplace convolution including the mapping property 
\begin{equation}
\label{GFI-map}
\I_{(\kappa)}:\, C_{-1}(0,+\infty)\, \rightarrow C_{-1}(0,+\infty),
\end{equation}
the commutativity law 
\begin{equation}
\label{GFI-com}
\I_{(\kappa_1)}\, \I_{(\kappa_2)} = \I_{(\kappa_2)}\, \I_{(\kappa_1)},\ \kappa_1,\, \kappa_2 \in \mathcal{L}_{1},
\end{equation}
and the index law
\begin{equation}
\label{GFI-index}
\I_{(\kappa_1)}\, \I_{(\kappa_2)} = \I_{(\kappa_1*\kappa_2)},\ \kappa_1,\, \kappa_2 \in \mathcal{L}_{1}.
\end{equation}

The first fundamental theorem of FC for the GFI \eqref{GFI}  and the GFDs \eqref{FDR-L} and \eqref{FD-C} of the Riemann-Liouville and the Caputo types, respectively, has been proved in \cite{Luc21a}.

\begin{theorem}[\cite{Luc21a}]
\label{t3}
Let $\kappa \in \mathcal{L}_{1}$ and $k$ be its associated Sonine kernel. 

Then,  the~GFD \eqref{FDR-L} is a left-inverse operator to the GFI \eqref{GFI} on the space $C_{-1}(0,+\infty)$:
\begin{equation}
\label{FTL}
(\D_{(k)}\, \I_{(\kappa)}\, f) (t) = f(t),\ f\in C_{-1}(0,+\infty),\ t>0,
\end{equation}
and the GFD \eqref{FD-C} is a left inverse operator to the GFI \eqref{GFI} on the space $C_{-1,(k)}^1(0,+\infty)$:
\begin{equation}
\label{FTC}
(\,_*\D_{(k)}\, \I_{(\kappa)}\, f) (t) = f(t),\ f\in C_{-1,(k)}^1(0,+\infty),\ t>0,
\end{equation}
where $C_{-1,(k)}^1(0,+\infty) := \{f:\ f(t)=(\I_{(k)}\, \phi)(t),\ \phi\in C_{-1}(0,+\infty)\}$.
\end{theorem}

In the rest of this section, we consider the compositions of the GFIs ($n$-fold GFIs) and construct the sequential GFDs. 

\begin{definition}[\cite{Luc21a}]
\label{d1}
Let $\kappa \in \mathcal{L}_{1}$. The $n$-fold GFI ($n \in \N$) is defined as a composition of $n$ GFIs with the kernel $\kappa$:
\begin{equation}
\label{GFIn}
(\I_{(\kappa)}^{<n>}\, f)(x) := (\underbrace{\I_{(\kappa)} \ldots \I_{(\kappa)}}_{n\ \mbox{times}}\, f)(x),\  x>0.
\end{equation}
\end{definition}

By employing the index law \eqref{GFI-index}, we can represent the $n$-fold GFI \eqref{GFIn} as a GFI with the kernel $\kappa^{<n>}$:
\begin{equation}
\label{GFIn-1}
(\I_{(\kappa)}^{<n>}\, f)(x) = (\kappa^{<n>}\, *\, f)(x) = (\I_{(\kappa)^{<n>}}\, f)(x),\  x>0.
\end{equation}

It is worth mentioning that the kernel $\kappa^{<n>},\ n\in \N$ belongs to the space 
$C_{-1}(0,+\infty)$, but it is not always a Sonine kernel. 

\begin{definition}
\label{d2}
Let $\kappa \in \mathcal{L}_{1}$ and $k$ be its associated Sonine kernel. The $n$-fold sequential GFDs in the Riemann-Liouville and the Caputo senses, respectively, are defined as follows:
\begin{equation}
\label{GFDLn}
(\D_{(k)}^{<n>}\, f)(x) := (\underbrace{\D_{(k)} \ldots \D_{(k)}}_{n\ \mbox{times}}\, f)(x),\  x>0,
\end{equation}
\begin{equation}
\label{GFDLn-C}
(\,_*\D_{(k)}^{<n>}\, f)(x) := (\underbrace{\,_*\D_{(k)} \ldots \,_*\D_{(k)}}_{n\ \mbox{times}}\, f)(x),\  x>0.
\end{equation}
\end{definition}

Please note that in \cite{Luc21b,Luc21a}, the $n$-fold  GFDs ($n \in \N$) were defined in a different form:
\begin{equation}
\label{GFDLn-1}
(\D_{(k)}^n\, f)(x) := \frac{d^n}{dx^n} ( k^{<n>} * f)(x),\ x>0,
\end{equation}
\begin{equation}
\label{GFDLn-1_C}
(\,_*\D_{(k)}^n\, f)(x) := ( k^{<n>} * f^{(n)})(x),\ x>0.
\end{equation}

The $n$-fold sequential GFDs \eqref{GFDLn} and \eqref{GFDLn-C} are generalizations of the Riemann-Liouville and the Caputo sequential fractional derivatives to the case of the Sonine kernels from $\mathcal{L}_{1}$.

Repeatedly applying the first fundamental theorem of FC for the GFI \eqref{GFI}  and the GFDs \eqref{FDR-L} and \eqref{FD-C} of the Riemann-Liouville and the Caputo type, respectively (Theorem \ref{t3}), we arrive at the following result:

\begin{theorem}[First Fundamental Theorem of FC for the $n$-fold sequential GFDs]
\label{t3-n}
Let $\kappa \in \mathcal{L}_{1}$ and $k$ be its associated Sonine kernel. 

Then,  the $n$-fold sequential GFD \eqref{GFDLn} in the Riemann-Liouville sense  is a left inverse operator to the $n$-fold  GFI \eqref{GFIn} on the space $C_{-1}(0,+\infty)$: 
\begin{equation}
\label{FTLn}
(\D_{(k)}^{<n>}\, \I_{(\kappa)}^{<n>}\, f) (x) = f(x),\ f\in C_{-1}(0,+\infty),\ x>0
\end{equation}
and the $n$-fold sequential GFD \eqref{GFDLn-C} in the Caputo sense  is a left inverse operator to the $n$-fold  GFI \eqref{GFIn} on the space $C_{-1,(k)}^n(0,+\infty)$: 
\begin{equation}
\label{FTLn-C}
(\,_*\D_{(k)}^{<n>}\, \I_{(\kappa)}^{<n>}\, f) (x) = f(x),\ f\in C_{-1,(k)}^n(0,+\infty),\ x>0,
\end{equation}
where $C_{-1,(k)}^n(0,+\infty) := \{f:\ f(x)=(\I_{(k)}^{<n>}\, \phi)(x),\ \phi\in C_{-1}(0,+\infty)\}$.
\end{theorem}

\section{Generalized convolution Taylor series}

For a Sonine kernel $\kappa \in \mathcal{L}_{1}$,  a convolution series in form \eqref{I-0-1} was introduced in \cite{Luc21c} in the framework of an operational calculus for the GFDs of the Caputo type with the Sonine kernels. It is worth mentioning that a part of the results regarding convolution series that were presented in \cite{Luc21c} is valid for any function $\kappa \in C_{-1}(0,+\infty)$ (that is not necessarily a Sonine kernel). In particular, this applies to Theorem 4.4 from \cite{Luc21c} that we now formulated and prove for a larger class of the kernels and in a slightly modified form.

\begin{theorem}
\label{t11}
Let a function $\kappa \in C_{-1}(0,+\infty)$  be represented in the form
\begin{equation}
\label{rep}
\kappa(x) = h_{p}(x)\kappa_1(x),\ x>0,\ p>0,\ \kappa_1\in C[0,+\infty)
\end{equation} 
and the convergence radius of the power  series
\begin{equation}
\label{ser}
\Sigma(z) = \sum^{+\infty }_{j=0}a_{j}\, z^j,\ a_{j}\in \C,\ z\in \C
\end{equation}
be non-zero. Then the convolution series 
\begin{equation}
\label{conser}
\Sigma_\kappa(x) = \sum^{+\infty }_{j=0}a_{j}\, \kappa^{<j+1>}(x)
\end{equation}
is convergent for all $x>0$ and defines a function  from the space $C_{-1}(0,+\infty)$. 
Moreover, the series 
\begin{equation}
\label{conser_p}
x^{1-\alpha}\, \Sigma_\kappa(x) = \sum^{+\infty }_{j=0}a_{j}\, x^{1-\alpha}\, \kappa^{<j+1>}(x),\ \ \alpha = \min\{p,\, 1\}
\end{equation}
is uniformly convergent for $x\in [0,\, X]$ for any $X>0$.
\end{theorem}

\begin{proof}

First we mention that in the case  $p\ge 1$ the function $\kappa$ is continuous on $[0,+\infty)$. Then the representation \eqref{rep} with $\kappa_1(x) = \kappa(x)$ and $p=1$ is valid. Thus, without any loss of generality, the representation \eqref{rep} holds valid with the parameter $p$ restricted to the interval $(0,1]$.

Now we introduce an arbitrary but fixed interval $[0,\, X]$ with $X > 0$. The function $\kappa_1(x) = \Gamma(p) x^{1-p}\kappa(x)$ from the representation \eqref{rep} is continuous on $[0,+\infty)$ and the estimate
\begin{equation}
\label{est_1}
\exists M_X >0:\ |\kappa_1(x)| = |\Gamma(p) x^{1-p}\kappa(x)| \le M_X,\ x \in [0,\, X]
\end{equation}
holds valid.

We proceed with derivation of a suitable estimate for the convolution powers $\kappa^{<j+1>},\ j\ge 1$ on the interval $(0,X]$. For $j=1$, we get
$$
|\kappa^{<2>}(x)|= |(\kappa*\kappa)(x)|\le \int_0^x h_{p}(x-\xi)|\kappa_1(x-\xi)|h_p(\xi)|\kappa_1(\xi)|\, d\xi \le
$$
$$
M_X^2  (h_{p}*h_{p})(x) = M_X^2\, h_{2p}(x),\ 0< x\le X.
$$
The same arguments easily lead to  the inequalities
\begin{equation}
\label{est_3}
|\kappa^{<j+1>}(x)|\le  M_X^{j+1}\, h_{(j+1)p}(x) = M_X^{j+1}\, \frac{x^{(j+1)p-1}}{\Gamma((j+1)p)},\ 0< x\le X,\ j\in \N_0 
\end{equation}
that can be rewritten as follows:
\begin{equation}
\label{est_3_3}
|x^{1-p}\, \kappa^{<j+1>}(x)|\le M_X^{j+1}\, \frac{x^{jp}}{\Gamma((j+1)p)},\ 0\le x\le X,\ j\in \N_0.
\end{equation}

For the further estimates, we choose and fix any point $z_0\not = 0$ from the convergence interval of the power series \eqref{ser}. Because the series  is absolutely convergent at the point $z_0$, the following inequalities hold true:
\begin{equation}
\label{est_2}
\exists M >0:\ |a_{j}\, z_0^j| \le M\ \forall j\in \N_0 \ \Rightarrow \
|a_{j}| \le \frac{M}{|z_0|^j}\ \forall j\in \N_0.
\end{equation}

Combining the inequalities \eqref{est_3_3} and \eqref{est_2}, we arrive  at  the following estimate:
\begin{equation}
\label{est_4}
|x^{1-p}\,a_j\, \kappa^{<j+1>}(x)|\le \frac{M}{|z_0|^j} M_X^{j+1}\, \frac{x^{jp}}{\Gamma((j+1)p)} \le
M\, M_X\, \frac{\left(\frac{M_X\, X^p}{|z_0|}\right)^j}{\Gamma((j+1)p)},\ j\in \N_0
\end{equation}
that is valid for all $x\in [0,\, X]$.
The number series
$$
\sum_{j=0}^{+\infty}M\, M_X\, \frac{\left(\frac{M_X\, X^p}{|z_0|}\right)^j}{\Gamma((j+1)p)}
$$
is absolutely convergent because of the Stirling asymptotic formula
$$
\Gamma(x+1) \sim \sqrt{2\pi x}\left( \frac{x}{e}\right)^x,\ x\to +\infty.
$$
This fact and the estimate \eqref{est_4} let us to apply the Weierstrass M-test that says that the series
\begin{equation}
\label{conser_1}
x^{1-p}\sum_{j=0}^{+\infty}a_j\kappa^{<j+1>}(x)
\end{equation}
is absolutely and uniformly convergent on the interval $[0,\, X]$. Because the functions $x^{1-p}a_j\kappa^{<j+1>}(x),\ j\in \N_0$ are all continuous on $[0,\, X]$ (see the inequality \eqref{est_3} and remember that $p\in (0,1]$), the uniform limit theorem ensures that the series \eqref{conser_1} is a function continuous on the interval $[0,\, X]$. Because $X$ can be chosen arbitrary large, the convolution series \eqref{conser}
is convergent for all $x>0$ and defines a function from the space $C_{-1}(0,+\infty)$.
\end{proof}

Now we proceed with analysis of the convolution series in form \eqref{conser} with the kernel functions $\kappa \in \mathcal{L}_{1}$. In what follows, we always assume that  the convergence radius of the power seres \eqref{ser} is non-zero. As proved in Theorem \eqref{t11}, the  convolution series \eqref{conser} defines a function from the space $C_{-1}(0,+\infty)$ that we denote by $f$:
\begin{equation}
\label{conser1}
f(x) = \sum^{+\infty }_{j=0}a_{j}\, \kappa^{<j+1>}(x).
\end{equation}

The problem that we now deal with is to determine the coefficients $a_j,\ j \in \N_0$ in terms of the function $f$. The series at the right-hand side of  \eqref{conser1} is uniformly convergent on any closed interval $[\delta,\, X],\ 0<\delta <X$ (see Theorem \ref{t11} and its proof) and thus for any $x>0$ we can apply the GFI \eqref{GFI} with the kernel $k$ (that is the Sonine kernel associated to $\kappa$) to this series term by term:
\begin{equation}
\label{a0} 
(\I_{(k)}\, f)(x) = (k\, * f)(x) = \left( k\, *\, \sum^{+\infty }_{j=0}a_{j}\, \kappa^{<j+1>}\right) (x) = 
\sum^{+\infty }_{j=0}a_{j}\, \left(k\, *\, \kappa^{<j+1>}\right) (x).
\end{equation}
Due to the Sonine condition \eqref{3-2}, the last formula can be represented in the form
\begin{equation}
\label{a0_1} 
(\I_{(k)}\, f)(x) = 
 a_0 + \left( \{1\} \, * \, \sum^{+\infty }_{j=0}a_{j+1}\, \kappa^{<j+1>}\right) (x) = a_0 + \left( \{1\} \, * \, f_1 \right) (x).
\end{equation}
According to Theorem \ref{t11}, the inclusion $f_1 \in C_{-1}(0,+\infty)$ holds valid.  As have been shown in \cite{LucGor99}, the definite integral of a function from $C_{-1}(0,+\infty)$ is a continuous function on the whole interval $[0,\, +\infty)$ that takes the value zero at the point zero:
\begin{equation}
\label{a0_2} 
\left( \{1\}\, * \, f_1 \right) (x) = (I_{0+}^1\, f_1)(x) \in C[0,\ +\infty),\ \ (I_{0+}^1\, f_1)(0) = 0.
\end{equation}
Substituting the point $x=0$ into the equation \eqref{a0_1}, we arrive at the formula
\begin{equation}
\label{a0_3} 
a_0 = (\I_{(k)}\, f)(0)
\end{equation}
for the coefficient $a_0$ of the convolution series \eqref{conser1}.

To proceed with the next coefficient, we first differentiate the representation \eqref{a0_1} and arrive at the formula
\begin{equation}
\label{a1_1} 
\frac{d}{dx}(\I_{(k)}\, f)(x) = 
\sum^{+\infty }_{j=0}a_{j+1}\, \kappa^{<j+1>}(x).
\end{equation}
The convolution series at the right-hand side of \eqref{a1_1} corresponds to the power series with the same convergence radius as the series \eqref{ser} and thus we can apply exactly same arguments as before to determine the coefficient $a_1$:
\begin{equation}
\label{a1_3} 
a_1 = \left(\I_{(k)}\, \frac{d}{dx}(\I_{(k)}\, f) \right)(0) = \left(\I_{(k)}\, \D_{(k)}\, f \right)(0).
\end{equation}

Repeating the same reasoning as for derivation of the coefficient $a_1$ again and again, we  get the formula
\begin{equation}
\label{an_3} 
a_j = \left(\I_{(k)}\, \D_{(k)}^{<j>}\, f \right)(0),\ j=2,3,\dots,
\end{equation}
where $\D_{(k)}^{<j>}$ stands for the $j$-fold sequential GFD  in the Riemann-Liouville sense  defined by 
\eqref{GFDLn}. Evidently, the formula \eqref{a0_3} is a particular case of the formula \eqref{an_3} with $j=0$. Summarizing the arguments presented above, we get a proof of the following theorem:

\begin{theorem}
\label{tgcTs}

Any function $f$ in form of a convolution series \eqref{conser1} with the Sonine kernel $\kappa \in \mathcal{L}_{1}$ can be represented as the following generalized convolution Taylor series:
\begin{equation}
\label{gcTs} 
f(x) = \sum^{+\infty }_{j=0}a_{j}\, \kappa^{<j+1>}(x),\ \ a_j = \left(\I_{(k)}\, \D_{(k)}^{<j>}\, f \right)(0),
\end{equation}
where $\I_{(k)}$ is the GFI \eqref{GFI}, $k$ is the Sonine kernel associated to the kernel $\kappa$ and $\D_{(k)}^{<j>}$ is the $j$-fold sequential GFD  in the Riemann-Liouville sense  defined by 
\eqref{GFDLn}.
\end{theorem}

\begin{example}
\label{ex1}
Let us illustrate the statement of Theorem \ref{tgcTs} on an example and consider the Sonine kernel $\kappa(x) = h_{\alpha}(x),\ 0 < \alpha <1$ with the associated kernel $k(x) = h_{1-\alpha}(x)$. Evidently, the inclusion $\kappa \in \mathcal{L}_{1}$ is valid and we can apply the result of Theorem \ref{tgcTs}. For the kernel $k(x) = h_{1-\alpha}(x),\ 0 < \alpha <1$, the GFI $\I_{(k)}$ is reduced to the Riemann-Lioville fractional integral of order $1-\alpha$ and the sequential GFD $\D_{(k)}^{<j>}$ is the well-known sequential Riemann-Lioville fractional derivative. As already mentioned, the convolution power $\kappa^{<j+1>}$ can be determined in explicit form:
\begin{equation}
\label{ex1_1}
\kappa^{<j+1>}(x) = h_{\alpha}^{<j+1>}(x) = h_{\alpha(j+1)}(x) = \frac{x^{\alpha(j+1)-1}}{\Gamma(\alpha(j+1))},\ x>0.
\end{equation}
The generalized convolution Taylor series \eqref{gcTs} takes then the following form:
\begin{equation}
\label{gcTs_ex1} 
f(x) = x^{\alpha-1}\, \sum^{+\infty }_{j=0}a_{j}\, \frac{x^{\alpha\, j}}{\Gamma(\alpha\, j +\alpha)},\ \ a_j = \left(I_{0+}^{1-\alpha}\, \left(D_{0+}^{\alpha}\right)^{<j>}\, f \right)(0),
\end{equation}
where $I_{0+}^{1-\alpha}$ is the Riemann-Lioville fractional integral \eqref{RLI} of order $1-\alpha$ and $\left(D_{0+}^{\alpha}\right)^{<j>}$ is the sequential Riemann-Lioville fractional derivative in form of a composition of $j$ Riemann-Lioville fractional derivatives \eqref{RLD} of order $\alpha$. The formula \eqref{gcTs_ex1} can be interpreted as a representation of a function $f$ from the space $C_{-1}(0,+\infty)$ in form of a power series with non-integer exponents. For other forms of such representations see e.g. \cite{Tru} and the references therein. It is worth mentioning that the representation \eqref{gcTs_ex1} is a particular case of a more general formula derived in  \cite{D58} for the case of the Dzherbashyan-Nersesyan fractional derivative. The general convolution Taylor series \eqref{gcTs} provides a far reaching generalizations of the results mentioned in this example.
\end{example}

\begin{example}
\label{ex2}
It is very instructive to look at the limiting case of the formula \eqref{gcTs_ex1} as $\alpha \to 1$. As already mentioned, the function $h_1(x) = \{1\}$ is not a Sonine kernel and thus in this case we cannot directly apply the theory presented above. However, we can use the relation \eqref{3-1-1} in  sense of generalized functions and derive the following conventional Taylor  series:
\begin{equation}
\label{gcTs_ex2} 
f(x) =  \sum^{+\infty }_{j=0}a_{j}\, \frac{x^{j}}{j!},\ \ a_j = \left(I_{0+}^{0}\, \left(D_{0+}^1\right)^{<j>}\, f \right)(0).
\end{equation}
Let us note that the coefficients $a_j$ coincide with those known in the theory of the conventional Taylor series (see the definitions \eqref{GFI_11}, \eqref{GFI_12}, \eqref{GFD_12_1}, and \eqref{GFD_12_2}):
$$
a_j = \left(I_{0+}^{0}\, \left(D_{0+}^1\right)^{<j>} f \right)(0)= (\delta\, * \, f^{(j)})(0) = f^{(j)}(0),\ j\in \N_0,
$$
$\delta$ being the Dirac $\delta$-function.
\end{example}

It is worth mentioning that in \cite{Luc21c} some important convolution series were introduced and employed for analytical treatment of the initial-value problems for the single- and multi-term fractional differential equations with the GFDs. 

If we start with the geometric series
\begin{equation}
\label{geom}
\Sigma(z) = \sum_{j=0}^{+\infty} \lambda^{j}z^{j},\ \lambda \in \C,\ z\in \C
\end{equation}
that for $\lambda \not =0$ has the convergence radius $r = 1/|\lambda|$, then  Theorem \ref{t11} ensures that the convolution series ($\kappa \in \mathcal{L}_{1}$)
\begin{equation}
\label{l}
l_{\kappa,\lambda}(x) = \sum_{j=0}^{+\infty} \lambda^{j}\kappa^{<j+1>}(x),\ \lambda \in \C
\end{equation}
is convergent for all $x>0$ and defines a function from the space $C_{-1}(0,+\infty)$.

For the kernel function $\kappa=\{ 1\}$, we immediately get the formula $\kappa^{<j+1>}(x) = \{ 1\}^{<j+1>}(x) = h_{j+1}(x)$. Then the convolution series \eqref{l} becomes a familiar power series for the exponential function:
\begin{equation}
\label{l-Mic}
l_{\kappa,\lambda}(x) = \sum_{j=0}^{+\infty} \lambda^{j}h_{j+1}(x) =
\sum_{j=0}^{+\infty} \frac{(\lambda\, x)^j}{j!} = e^{\lambda\, x}.
\end{equation}

In the case of the kernel $\kappa(x) = h_{\alpha}(x)$ of the Riemann-Liouville fractional integral,  the formula  $\kappa^{<j+1>}(x) = h_{\alpha}^{<j+1>}(x) = h_{(j+1)\alpha}(x)$ is valid and the convolution series \eqref{l} takes the form
\begin{equation}
\label{l-Cap}
l_{\kappa,\lambda}(x) = \sum_{j=0}^{+\infty} \lambda^{j}h_{(j+1)\alpha}(x) =
x^{\alpha-1}\sum_{j=0}^{+\infty} \frac{\lambda^j\, x^{j\alpha}}{\Gamma(j\alpha+\alpha)} = x^{\alpha -1}E_{\alpha,\alpha}(\lambda\, x^{\alpha}),
\end{equation}
where the two-parameters Mittag-Leffler function $E_{\alpha,\alpha}$ is defined by \eqref{ML}.

Another interesting case is the kernel $\kappa(x) = h_{1-\beta+\alpha}(x)\, +\, h_{1-\beta}(x),\ 0<\alpha < \beta <1$ (see the formula \eqref{3-8}). As shown in \cite{Luc21c}, the convolution series \eqref{l} takes in this case the following form:
$$
l_{\kappa,\lambda}(x) = 
\frac{1}{\lambda x} \sum_{j=0}^{+\infty}\sum_{l_1+l_2 =j} \frac{j!}{l_1!l_2!}\frac{(\lambda x^{1-\beta+\alpha})^{l_1}(\lambda x^{1-\beta})^{l_2}}{\Gamma(l_1(1-\beta+\alpha)+l_2(1-\beta))} =
$$
$$
\frac{1}{\lambda x}  E_{(1-\beta,1-\beta+\alpha),0}(\lambda x^{1-\beta},\lambda x^{1-\beta+\alpha}),
$$
where $E_{(1-\beta,1-\beta+\alpha),0}$ is a particular case of the multinomial Mittag-Leffler function
\begin{equation}
\label{MLm}
E_{(\alpha_1,\ldots,\alpha_m),\beta}(z_1,\ldots,z_m):=
\sum_{j=0}^{+\infty} \sum_{l_1+\cdots +l_m =j}\frac{j!}{l_1!\times\cdot\times
l_m!}\frac{\prod_{i=1}^m
z_i^{l_i}}{\Gamma(\beta+\sum_{i=1}^m \alpha_i l_i)}
\end{equation}
introduced for the first time in \cite{Luc93} (see also \cite{HadLuc}).

In \cite{Luc21c}, the convolution series of type \eqref{l} and their generalizations were employed for derivation of analytical solutions to the initial-value problems for the fractional differential equation with the GFDs of Caputo type. 

\section{Generalized convolution Taylor formula}

In this section, we derive two different forms of the generalized convolution Taylor formula for  representation of the functions from a certain space in form of the convolution polynomials with a remainder in terms of the GFIs and GFDs. 

We start with the case of the GFD defined in the Riemann-Liouville sense and first prove the following result formulated for the functions from the space $C_{-1,(k)}^{(1)}(0,+\infty) = \{ f\in C_{-1}(0,+\infty):\ (\D_{(k)}\, f)\in C_{-1}(0,+\infty)\}$ (please note that this space does not coincide with the space $C_{-1,(k)}^{1}(0,+\infty)$ introduced in the previous section and the inclusion $C_{-1,(k)}^{1}(0,+\infty) \subset C_{-1,(k)}^{(1)}(0,+\infty)$ is valid). 

\begin{theorem}[Second Fundamental Theorem of FC for the GFD in the R-L sense]
\label{tgcTf}
Let $\kappa \in \mathcal{L}_{1}$ and $k$ be its associated Sonine kernel. 
For a function $f\in C_{-1,(k)}^{(1)}(0,+\infty)$, the formula
\begin{equation}
\label{sFTL}
(\I_{(\kappa)}\, \D_{(k)}\, f) (x) =  f(x) - (k\, *\, f)(0)\kappa(x) = f(x) - (\I_{(k)}\, f)(0)\kappa(x),\ x>0
\end{equation}
holds valid.
\end{theorem}
\begin{proof}
First we determine the kernel of the GFD $\D_{(k)}$ on the space $C_{-1,(k)}^{(1)}(0,+\infty)$:
$$
(\D_{(k)}\, f) (x) = \frac{d}{dx} ((\I_{(k)}\, f) (x) = 0 \ \Leftrightarrow \ (\I_{(k)}\, f) (x) = C,\ C\in \R \ \Leftrightarrow 
$$
$$
(k\, *\,  f) (x) = C \ \Leftrightarrow \ (\kappa \, *\, (k\, *\,  f)) (x) = (\kappa \, *\, \{C\})(x) = C  (\kappa \, *\, \{1\})(x) \ \Leftrightarrow \ 
$$
$$
(\{1\} \, *\, f)(x) = (\{1\} \, *\, C\kappa )(x) \ \Leftrightarrow \ f(x) = C\, \kappa(x).
$$
Thus, the kernel of $\D_{(k)}$ is as follows:
\begin{equation}
\label{kern}
\mbox{Ker}\, \D_{(k)} = \{C\, \kappa(x):\, C\in \R\}.
\end{equation}
Let us now introduce an auxiliary function
\begin{equation}
\label{F}
F(x) := (\I_{(\kappa)}\, \D_{(k)}\, f) (x).
\end{equation}
Because the GFD \eqref{FDR-L} is a left inverse operator to the GFI \eqref{GFI} on the space $C_{-1}(0,+\infty)$ (see  
the formula \eqref{FTL}), we get the relation
$$
(\D_{(k)}\, F)(x) = (\D_{(k)}\, \I_{(\kappa)}\, \D_{(k)}\, f) (x) = (\D_{(k)}\, f) (x)
$$
and thus the function $\psi(x) = F(x)-f(x)$ belongs to the kernel of the GFD $\D_{(k)}$, i.e.,
\begin{equation}
\label{psi}
\psi(x) = F(x)-f(x) = C\, \kappa(x).
\end{equation}
To determine the constant $C$, we act on the last relation  with the GFI $\I_{(k)}$ and get the formula
$$
(\I_{(k)}\, (F-f))(x) = (k\, *\, C\kappa)(x) = C \{1\} = C.
$$
Otherwise,
$$
(\I_{(k)}\, F)(x)= (\I_{(k)}\,(\I_{(\kappa)}\, \D_{(k)}\, f)) (x) = ((\I_{(k)}\, \I_{(\kappa)})\, \D_{(k)}\, f) (x) = 
(\{1\}\, *\,  \D_{(k)}\, f) (x).
$$
Combining the last two relations, we arrive at the formula
\begin{equation}
\label{kern1}
(\{1\}\, *\,  \D_{(k)}\, f) (x) - (\I_{(k)}\, f)(x) = C. 
\end{equation}
For a function $f\in C_{-1,(k)}^{(1)}(0,+\infty)$, the inclusion $\D_{(k)}\, f \in C_{-1}(0,+\infty)$ holds valid. Thus 
$\{1\}\, *\,  \D_{(k)}\, f \in C[0,+\infty)$ and $(\{1\}\, *\,  \D_{(k)}\, f) (0) = 0$. Substituting $x=0$ into the formula \eqref{kern1} leads to the relation
$$
C = - (\I_{(k)}\, f)(0)
$$
that together with the formulas \eqref{F} and \eqref{psi} finalizes  the proof of the theorem. 
\end{proof}

For the Sonine kernel $\kappa(x)= h_{\alpha}(x),\ 0 <\alpha <1$, the representation \eqref{sFTL} is well-known (see e.g. \cite{Samko}):
\begin{equation}
\label{sFTL-a}
(I_{0+}^\alpha\, D_{0+}^\alpha\, f) (x) = f(x) - (I_{0+}^{1-\alpha}\, f)(0)\frac{x^{\alpha -1}}{\Gamma(\alpha)},\ x>0,
\end{equation}
where $I_{0+}^\alpha$ and $D_{0+}^\alpha$ are the Riemann-Liouville fractional integral and derivative, respectively. 

It is also worth mentioning that  in the case $\kappa(x)= h_{1}(x)$ the space of functions $C_{-1,(k)}^{(1)}(0,+\infty)$  corresponds to the space of continuously differentiable functions and the formula \eqref{sFTL} reads
$$  
\int_0^x \, f^\prime(\xi)\, d\xi = f(x) - f(0).
$$

Now we generalize Theorem \ref{tgcTf} to the case of the $n$-fold GFIs and the $n$-fold sequential GFDs in the Riemann-Liouville sense. This time, the result is formulated for the functions from the space $C_{-1,(k)}^{(n)}(0,+\infty) = \{ f\in C_{-1}(0,+\infty):\ (\D_{(k)}^{<j>}\, f)\in C_{-1}(0,+\infty),\ j=1,\dots,n\}$ that in the case $\kappa(x)= h_{1}(x)$ corresponds to the  space of $n$-times continuously differentiable functions. 

\begin{theorem}[Second Fundamental Theorem of FC for the sequential GFD in the R-L  sense]
\label{tgcTfn}
Let $\kappa \in \mathcal{L}_{1}$ and $k$ be its associated Sonine kernel. 
For a function $f\in C_{-1,(k)}^{(n)}(0,+\infty)$, the formula
\begin{equation}
\label{sFTLn}
(\I_{(\kappa)}^{<n>}\, \D_{(k)}^{<n>}\, f) (x) = f(x) - \sum_{j=0}^{n-1}\left( k\, * \, \D_{(k)}^{<j>}\, f\right)(0)\kappa^{<j+1>}(x) = 
\end{equation}
$$
f(x) - \sum_{j=0}^{n-1}\left( \I_{(k)}\, \D_{(k)}^{<j>}\, f\right)(0)\kappa^{<j+1>}(x),\ x>0
$$
holds valid, where $\I_{(\kappa)}^{<n>}$ is the $n$-fold GFI \eqref{GFIn} and $\D_{(k)}^{<n>}$ is the $n$-fold sequential GFD \eqref{GFDLn} in the Riemann-Liouville sense. 
\end{theorem}

\begin{proof}
To prove the formula \eqref{sFTLn}, we repeatedly employ Theorem \ref{tgcTf}. For $n=2$, we first get the representation
$$ 
(\I_{(\kappa)}^{<2>}\, \D_{(k)}^{<2>}\, f) (x) = (\I_{(\kappa)} \I_{(\kappa)}\, \D_{(k)} \, \D_{(k)}\, f) (x) = 
$$
$$
(\I_{(\kappa)} (\I_{(\kappa)}\, \D_{(k)} \, (\D_{(k)}f))) (x).
$$
Then we apply Theorem \eqref{tgcTf} to the inner composition   $\I_{(\kappa)}\, \D_{(k)}$ acting on the function $(\D_{(k)}\, f)$ and get the formula
$$
(\I_{(\kappa)}^{<2>}\, \D_{(k)}^{<2>}\, f) (x) = \left(\I_{(\kappa)}\, \left[ (\D_{(k)}\, f)(x) - \left( k\, *  \D_{(k)}\, f\right)(0)\, \kappa(x) \right]\right)(x) = 
$$
$$
(\I_{(\kappa)}\, \D_{(k)}\, f)(x) - \left( k\, *  \D_{(k)}\, f\right)(0)\, \kappa^{<2>}(x).
$$
Now we apply Theorem \eqref{tgcTf} once again, this time to the composition $\I_{(\kappa)}\, \D_{(k)}\, f$ at the right-hand side of the last formula and get the final result:
$$
(\I_{(\kappa)}^{<2>}\, \D_{(k)}^{<2>}\, f) (x) = f(x) - (k\, *\, f)(0)\kappa(x) - \left( k\, *  \D_{(k)}\, f\right)(0)\, \kappa^{<2>}(x).
$$
To proceed with the general case, we  employ the following recurrent formula:
$$
(\I_{(\kappa)}^{<n>}\, \D_{(k)}^{<n>}\, f) (x) = (\I_{(\kappa)}^{<n-1>} (\I_{(\kappa)}\, \D_{(k)} \, (\D_{(k)}^{<n-1>}f))) (x) =
$$
$$
\left(\I_{(\kappa)}^{<n-1>} \left[ (\D_{(k)}^{<n-1>}f)(x) - \left( k\, * \, \D_{(k)}^{<n-1>}f\right)(0)\kappa(x)\right]\right)(x),\ n=3,4,\dots.
$$
The representation \eqref{sFTLn} easily follows from this recurrent formula and the principle of the mathematical induction. 
\end{proof}

Now we are ready to formulate and prove one of the main results of this section in form of the following theorem:

\begin{theorem}[Generalized convolution Taylor formula for the GFD in the R-L sense]
\label{tgcTfn-1}
Let $\kappa \in \mathcal{L}_{1}$ and $k$ be its associated Sonine kernel. 
For a function $f\in C_{-1,(k)}^{(n)}(0,+\infty)$, the generalized convolution Taylor formula
\begin{equation}
\label{gcTf}
f(x) = \sum_{j=0}^{n-1}a_j\, \kappa^{<j+1>}(x) + r_n(x), \ x>0 
\end{equation}
holds valid. The coefficients $a_j$ are given by the expression
\begin{equation}
\label{aj}
a_j = \left( k\, * \, \D_{(k)}^{<j>}\, f\right)(0) =  \left( \I_{(k)}\, \D_{(k)}^{<j>}\, f\right)(0),\ j=0,1\dots,n-1
\end{equation}
and the remainder can be represented as follows:
\begin{equation}
\label{rem}
r_n(x) = 
(\I_{(\kappa)}^{<n>}\, \D_{(k)}^{<n>}\, f) (x) = (\D_{(k)}^{<n>}\, f) (\xi) \, (\{ 1\} \, *\, \kappa^{<n>})(x),\ 0< \xi \le x,
\end{equation}
where $\I_{(\kappa)}^{<n>}$ is the $n$-fold GFI \eqref{GFIn} and $\D_{(k)}^{<n>}$ is the $n$-fold sequential GFD \eqref{GFDLn} in the Riemann-Liouville sense. 
\end{theorem}
\begin{proof}
First we mention that the generalized convolution Taylor formula \eqref{gcTf} with the coefficients \eqref{aj} and the remainder in form $r_n(x) = (\I_{(\kappa)}^{<n>}\, \D_{(k)}^{<n>}\, f) (x)$ immediately follows from 
the second Fundamental Theorem of FC for the sequential GFD in the Riemann-Liouville sense (Theorem \ref{tgcTfn}). The second form of the remainder is obtained by application of the integral mean value theorem:
$$
(\I_{(\kappa)}^{<n>}\, \D_{(k)}^{<n>}\, f) (x) = \int_0^x \kappa^{<n>}(x-\xi)\, (\D_{(k)}^{<n>}\, f) (\xi)\, d\xi = 
$$
$$
(\D_{(k)}^{<n>}\, f) (\xi) \int_0^x \kappa^{<n>}(x-\xi)\, \, d\xi = \D_{(k)}^{<n>}\, f) (\xi) \, (\{ 1\} \, *\, \kappa^{<n>})(x),\ 0<\xi \le x.
$$
\end{proof}

In the case $\kappa(x) = h_\alpha(x),\ 0<\alpha <1$ (the kernel of the Riemann-Liouville fractional integral), the generalized convolution Taylor formula \eqref{gcTf} is reduced to the following known form (\cite{D58}):
\begin{equation}
\label{gcTf-1}
f(x) = x^{\alpha-1} \sum_{j=0}^{n-1}a_j\, \frac{x^{\alpha\, j}}{\Gamma(\alpha\, j +\alpha)} + r_n(x), \ x>0, 
\end{equation}
\begin{equation}
\label{aj-1}
a_j =  \left( I_{0+}^{1-\alpha}\, \left(D_{0+}^\alpha\right)^{<j>}\, f\right)(0),\ j=0,1\dots,n-1
\end{equation}
with the remainder 
\begin{equation}
\label{rem-1}
r_n(x) = 
\left(I_{0+}^{n\, \alpha}\, \left(D_{0+}^\alpha\right)^{<n>}\, f\right) (x) = \left(\left(D_{0+}^\alpha\right)^{<n>}\, f\right) (\xi)\frac{x^{\alpha\, n}}{\Gamma(\alpha\, n +\alpha)},\ 0< \xi \le x,
\end{equation}
where $I_{0+}^{\alpha}$ is the Riemann-Liouville fractional integral and $\left(D_{0+}^\alpha\right)^{<n>}$ is the $n$-fold sequential Riemann-Liouville fractional derivative. The conventional Taylor formula is obtained from the formula \eqref{gcTf-1} by letting $\alpha$ go to $1$. 

Now we proceed with the case of the GFD in the Caputo sense. This time, the results will be formulated for the functions from the space $C_{-1}^n(0,+\infty) = \{ f:\ f^{(n)}\in C_{-1}(0,+\infty)\},\ n\in \N$. This space of functions was introduced in \cite{LucGor99} and employed in \cite{Luc21b,Luc21a,Luc21c} for derivation of several results regarding the GFD in the Caputo sense. In particular, the following result was proved in \cite{Luc21a}:

\begin{theorem}[Second Fundamental Theorem of FC for the GFD in the Caputo sense]
\label{t4}
Let $\kappa \in \mathcal{L}_{1}$ and $k$ be its associated Sonine kernel. 

Then,  the~relation
\begin{equation}
\label{2FTC}
(\I_{(\kappa)}\,_*\D_{(k)}\, f) (x) = f(x)-f(0),\ x>0
\end{equation}
holds valid for the functions $f\in C_{-1}^1(0,+\infty)$.
\end{theorem}

For the Sonine kernel $\kappa(x)= h_{\alpha}(x),\ 0 <\alpha <1$ that generates the Caputo fractional derivative  $\,_*D_{0+}^\alpha$ defined by 
\eqref{CD}, the formula \eqref{2FTC} is well-known (see e.g. \cite{LucGor99}):
\begin{equation}
\label{2FTL-a}
(I_{0+}^\alpha\,_*D_{0+}^\alpha\, f) (x) = f(x) - f(0),\ x>0,
\end{equation}
where $I_{0+}^\alpha$ is the Riemann-Liouville fractional integral. 

As in the case of the GFD in the Riemann-Liouville sense,  we generalize Theorem \ref{t4} to the case of the $n$-fold GFI and the $n$-fold sequential GFD in the Caputo sense. 

\begin{theorem}[Second Fundamental Theorem of FC for the sequential GFD in the Caputo sense]
\label{tgcTfn-C}
Let $\kappa \in \mathcal{L}_{1}$ and $k$ be its associated Sonine kernel. 
For a function $f\in C_{-1}^{n}(0,+\infty)$, the formula
\begin{equation}
\label{sFTLn-C}
(\I_{(\kappa)}^{<n>}\,_*\D_{(k)}^{<n>}\, f) (x) = f(x) - f(0) - \sum_{j=1}^{n-1}\left(\,_*\D_{(k)}^{<j>}\, f\right)(0)\left( \{1\} \, *\, \kappa^{<j>}\right)(x)
\end{equation}
holds valid, where $\I_{(\kappa)}^{<n>}$ is the $n$-fold GFI \eqref{GFIn} and $\,_*\D_{(k)}^{<n>}$ is the $n$-fold sequential GFD \eqref{GFDLn-C}. 
\end{theorem}

\begin{proof}
The formula \eqref{sFTLn-C} immediately follows from Theorem \ref{t4}. Indeed, for $n=2$, we first get the representation
$$ 
(\I_{(\kappa)}^{<2>}\,_*\D_{(k)}^{<2>}\, f) (x) = (\I_{(\kappa)} \I_{(\kappa)}\,_*\D_{(k)}\,_*\D_{(k)}\, f) (x) = 
$$
$$
(\I_{(\kappa)} (\I_{(\kappa)}\,_*\D_{(k)} \, (\,_*\D_{(k)}f))) (x).
$$
Then we apply Theorem \ref{t4} and get the formula
$$
(\I_{(\kappa)}^{<2>}\,_*\D_{(k)}^{<2>}\, f) (x) = \left(\I_{(\kappa)}\, \left[ (\,_*\D_{(k)}\, f)(x) - \left(\,_*\D_{(k)}\, f\right)(0)\right]\right)(x) = 
$$
$$
(\I_{(\kappa)}\,_*\D_{(k)}\, f)(x) - ( \I_{(\kappa)} (\,_*\D_{(k)}\, f)(0))(x) =
$$
$$
f(x) - f(0) - (\,_*\D_{(k)}\, f)(0)(\{1\}\, *\, \kappa)(x).
$$
In the general case, the representation \eqref{sFTLn-C} immediately follows from the recurrent formula
$$
(\I_{(\kappa)}^{<n>}\,_*\D_{(k)}^{<n>}\, f) (x) = (\I_{(\kappa)}^{<n-1>} (\I_{(\kappa)}\,_*\D_{(k)} \, (\,_*\D_{(k)}^{<n-1>}f))) (x) =
$$
$$
\left(\I_{(\kappa)}^{<n-1>} \left[ (\,_*\D_{(k)}^{<n-1>}f)(x) - \left( \,_*\D_{(k)}^{<n-1>}f\right)(0)\right]\right)(x),\ n=3,4,\dots.
$$
and the principle of the mathematical induction.  
\end{proof}

The representation \eqref{sFTLn-C} can be rewritten in form of a generalized convolution Taylor formula.

\begin{theorem}[Generalized convolution Taylor formula for the GFD in the Caputo sense]
\label{tgcTfn-1-C}
Let $\kappa \in \mathcal{L}_{1}$ and $k$ be its associated Sonine kernel. 
For a function $f\in C_{-1}^{n}(0,+\infty)$, the generalized convolution Taylor formula
\begin{equation}
\label{gcTf-C}
f(x) = f(0) + \sum_{j=1}^{n-1}\left(\,_*\D_{(k)}^{<j>}\, f\right)(0)\, \left(\{1\}\, *\,  \kappa^{<j>}\right)(x) + R_n(x), \ x>0 
\end{equation}
holds valid with the remainder in the form
\begin{equation}
\label{rem-C}
R_n(x) = 
(\I_{(\kappa)}^{<n>}\,_*\D_{(k)}^{<n>}\, f) (x) =  (\,_*\D_{(k)}^{<n>}\, f) (\xi) \, (\{ 1\} \, *\, \kappa^{<n>})(x),\ 0< \xi \le x,
\end{equation}
where $\I_{(\kappa)}^{<n>}$ is the $n$-fold GFI \eqref{GFIn} and $\,_*\D_{(k)}^{<n>}$ is the $n$-fold sequential GFD \eqref{GFDLn-C} in the Caputo sense. 
\end{theorem}

The second form of the remainder is obtained by application of the integral mean value theorem:
$$
(\I_{(\kappa)}^{<n>}\,_*\D_{(k)}^{<n>}\, f) (x) = \int_0^x \kappa^{<n>}(x-\xi)\, (\,_*\D_{(k)}^{<n>}\, f) (\xi)\, d\xi = 
$$
$$
(\,_*\D_{(k)}^{<n>}\, f) (\xi) \int_0^x \kappa^{<n>}(x-\xi)\, \, d\xi = (\,_*\D_{(k)}^{<n>}\, f) (\xi) \, (\{ 1\} \, *\, \kappa^{<n>})(x),\ 0<\xi \le x.
$$

In the case $\kappa(x) = h_\alpha(x),\ 0<\alpha <1$ (the kernel of the Caputo fractional derivative), $\kappa^{<n>}(x) =h_{\alpha\, n}(x)$ and the generalized convolution Taylor formula \eqref{gcTf-C} takes the form:
\begin{equation}
\label{gcTf-1-C}
f(x) = f(0) +  \sum_{j=1}^{n-1}\left(\left(\,_*D_{0+}^\alpha\right)^{<j>}\, f\right)(0)\, \frac{x^{\alpha\, j}}{\Gamma(\alpha\, j +1)} + R_n(x), \ x>0 
\end{equation}
with the remainder 
\begin{equation}
\label{rem-1-C}
R_n(x) = 
\left(I_{0+}^{n\, \alpha}\, \left(\,_*D_{0+}^\alpha\right)^{<n>}\, f\right) (x),
\end{equation}
where $I_{0+}^{\alpha}$ is the Riemann-Liouville fractional integral and $\left(\,_*D_{0+}^\alpha\right)^{<n>}$ is the $n$-fold sequential Caputo fractional derivative. 

As we see, the generalized convolution Taylor formula involving the Riemann-Liouville fractional derivative (formula \eqref{gcTf-1}) and the one involving the Caputo fractional derivative (formula \eqref{gcTf-1-C}) are completely different. Whereas the generalized convolution Taylor polynomial at the right-hand side of the formula \eqref{gcTf-1} has an integrable singularity of a power function type at the origin, the Taylor polynomial at right-hand side of the formula \eqref{gcTf-1-C} is continuous at the point $x=0$. 

Letting $n$ go to $+\infty$ in the formula \eqref{gcTf-C} leads to the generalized convolution Taylor series for the GFD in the Caputo sense in the form
\begin{equation}
\label{GC}
f(x) = f(0) +  \sum_{j=1}^{+\infty}\left(\,_*\D_{(k)}^{<j>}\, f\right)(0)\,\left( \{1\}\, *\, \kappa^{<j>}\right)(x), \ x>0. 
\end{equation}
The representation \eqref{GC} holds valid under the condition that the convolution series at its right-hand side is a convergent one (this is the case, say, if the sequence $\left(\,_*\D_{(k)}^{<j>}\, f\right)(0),\ j\in \N$ is bounded).

The generalized convolution Taylor formulas and the generalized convolution Taylor series that were presented in this section can be applied among other things for derivation of analytical solutions to the fractional differential equations with the GFDs in the Riemann-Liouville sense. This type of equations was not yet treated in the FC literature and will be considered elsewhere.


\end{document}